\documentclass{article}
\usepackage{amssymb,latexsym,amsmath,amsthm,mathrsfs}
\usepackage{graphicx}
\newcommand{\comment}[1]{}

\newcommand{\combin}[2]{\bigg(\frac{#1}{#2}\bigg)}

\DeclareMathOperator{\Cos}{cos.}
\DeclareMathOperator{\Sin}{sin.}

\begin{document}
\title{General observations on series whose terms proceed as the sines
and cosines of multiples of angles\footnote{Presented to the St. Petersburg Academy
on March 6, 1777. Originally published as
{\em Observationes generales circa series, quarum termini secundum sinus vel
cosinus angulorum multiplorum progrediuntur},
Nova acta academiae scientiarum Petropolitinae \textbf{7} (1789),
87--98.
E655 in the Enestr{\"o}m index.
Translated from the Latin by Jordan Bell,
Department of Mathematics, University of Toronto, Toronto, Ontario, Canada.
Email: jordan.bell@utoronto.ca}}
\author{Leonhard Euler}
\date{}
\maketitle

\S 1. If the summation of the series
\[
A+Bx+Cxx+Dx^3+\textrm{etc.}
\]
were known, so that whatever value is taken for the letter $x$, the sum
can be assigned, then also the sum of the series
\[
A+B\Cos\varphi+C\Cos 2\varphi+D\Cos 3\varphi+\textrm{etc.}
\]
and
\[
B\Sin \varphi+C\Sin 2\varphi+D\Sin 3\varphi+E\Sin 4\varphi+\textrm{etc.}
\]
can always be exhibited. If the sum of the first series is expressed
as a certain function of $x$, which we shall designate by the character
$\Delta:x$, so that it is
\[
\Delta:x=A+Bx+Cxx+Dx^3+\textrm{etc.},
\]
if we write in place of $x$
\[
\Cos\varphi+\surd -1 \cdot \Sin \varphi
\]
or 
\[
\Cos\varphi-\surd -1\cdot \Sin\varphi,
\]
the sum of the resulting series will be
\[
2A+2B\Cos\varphi+2C\Cos 2\varphi+2D\Cos 3\varphi+2E\Cos 4\varphi+\textrm{etc.}
\]
of which the sum will thus be
\[
\Delta:(\Cos\varphi+\surd -1 \cdot \Sin\varphi)
+\Delta:(\Cos\varphi-\surd -1\cdot \Sin\varphi);
\]
but on the other hand let us subtract the second from the first, which
will produce this series
\[
2B\surd -1 \cdot\Sin\varphi+2C\surd -1\cdot \Sin 2\varphi
+2D\surd -1\cdot \Sin 3\varphi
+2E\surd -1\cdot \Sin 4\varphi+\textrm{etc.}
\]
whose sum will therefore be
\[
\Delta:(\Cos \varphi+\surd -1\cdot \Sin \varphi)
-\Delta:(\Cos \varphi-\surd -1\cdot
\Sin \varphi).
\]

\S 2. To make these expressions more convenient to work with,
let us set
for the sake of brevity
\[
\Cos \varphi+\surd -1\cdot\Sin\varphi=p,
\quad \Cos\varphi-\surd -1\cdot\Sin\varphi=q
\]
and it will be, as is commonly known,
\[
pq=1
\]
and so
\[
q=\frac{1}{p};
\]
then indeed it will be
\[
\Cos\varphi=\frac{p+q}{2}; \quad \Cos 2\varphi=\frac{pp+qq}{2}; \quad
\Cos 3\varphi=\frac{p^3+q^3}{2}; \quad
\Cos 4\varphi=\frac{p^4+q^4}{2} \quad \textrm{etc.}
\]
Also indeed for sines it will be obtained that
\[
\Sin\varphi=\frac{p-q}{2\surd -1}; \quad \Sin 2\varphi=\frac{pp-qq}{2\surd -1};
\quad \Sin 3\varphi=\frac{p^3-q^3}{2\surd -1}; \quad
\Sin 4\varphi=\frac{p^4-q^4}{2\surd -1} \quad \textrm{etc.}
\]
Having established this, these two summations are obtained
\[
A\Cos 0\varphi+B\Cos \varphi+C\Cos 2\varphi+D\Cos 3\varphi
+\textrm{etc.}=
\frac{\Delta:p+\Delta:q}{2};
\]
and
\[
A\Sin 0\varphi+B\Sin \varphi+C\Sin 2\varphi+D\Sin 3\varphi+\textrm{etc.}
=\frac{\Delta:p-\Delta:q}{2\surd -1}.
\]

\S 3. Let us now take for the principal series any power of the binomial
expansion, which is
\[
(1+x)^n=1+\frac{n}{1}x+\frac{n(n-1)}{1\cdot 2}xx+\frac{n(n-1)(n-2)}{1\cdot 2\cdot 3}x^3+\textrm{etc.},
\]
so that in this case it would be
\[
\Delta:x=(1+x)^n;
\]
then indeed to shorten this expression, we shall designate all the
coefficients by these characters,
which we have already used several times,
\[
\combin{n}{0}, \combin{n}{1}, \combin{n}{2}, \combin{n}{3}, \combin{n}{4} \textrm{etc.}
\]
so that it is
\begin{eqnarray*}
\combin{n}{0}&=&1,\\
\combin{n}{1}&=&n,\\
\combin{n}{2}&=&\frac{n}{1}\cdot \frac{n-1}{2},\\
\combin{n}{3}&=&\frac{n}{1}\cdot \frac{n-1}{2}\cdot \frac{n-2}{3}\\
&&\textrm{etc.},
\end{eqnarray*}
where it is useful to observe in general that
\[
\combin{n}{i}=\combin{n}{n-i}
\]
and hence
\[
\combin{n}{n}=\combin{n}{0}=1.
\]
As well it is indeed evident that whenever $i$ is a negative number or
a positive number greater than $n$, it will always be
\[
\combin{n}{i}=0,
\]
if indeed $n$ is an integral number. Having made these observations, we will
thus have the principal summation 
\[
(1+x)^n=\combin{n}{0}+\combin{n}{1}x+\combin{n}{2}x^2+\combin{n}{3}x^3
+\combin{n}{4}x^4+\textrm{etc.},
\]
from which by the precept just given we will derive the following two summations
\[
\combin{n}{0}\Cos 0\varphi+\combin{n}{1}\Cos\varphi+\combin{n}{2}\Cos 2\varphi+
\combin{n}{3}\Cos 3\varphi+\textrm{etc.}=\frac{(1+p)^n+(1+q)^n}{2}
\]
and
\[
\combin{n}{0}\Sin 0\varphi+\combin{n}{1}\Sin \varphi+\combin{n}{2}\Sin 2\varphi
+\combin{n}{3} \Sin 3\varphi+\textrm{etc.}=\frac{(1+p)^n-(1+q)^n}{2\surd -1}.
\]
For any case where the formulae for $p$ and $q$ are assumed imaginary,
the above formulae will always return real values, which we shall
show in the following two problems.

\begin{center}
{\Large Problem 1}
\end{center}

{\em Given this series of cosines
\[
1+\frac{n}{1}\Cos \varphi+\frac{n}{1}\cdot \frac{n-1}{2}\Cos 2\varphi
+\frac{n}{1}\cdot \frac{n-1}{2}\cdot \frac{n-2}{3}\Cos 3\varphi+\textrm{etc.}
=s,
\]
which with the characters that have been established will be
\[
s=\combin{n}{0}\Cos 0\varphi+\combin{n}{1}\Cos \varphi
+\combin{n}{2}\Cos 2\varphi+\combin{n}{3}\Cos 3\varphi+\textrm{etc.},
\]
to express its sum with reals.}

\begin{center}
{\Large Solution}
\end{center}

\S 4. As we observed,
\[
s=\frac{(1+p)^n+(1+q)^n}{2},
\]
with it being
\[
p=\Cos \varphi+\surd -1\cdot \Sin \varphi \quad \textrm{and}
\quad q=\Cos \varphi-\surd -1\cdot \Sin \varphi,
\]
the whole matter is reduced to this, that the expression for exhibiting $s$
is freed from imaginaries; for it is evident
here that the imaginaries mutually eliminate each other, which will
appear by summing; 
thus we will inquire into another resolution such that 
the imaginaries can be removed without
doing the expansion, 
which will be done in the following way.

\S 5. Since it is $pq=1$, the formula $1+p$ can thus be expressed as such
\[
1+p=(\surd p+\surd q)\surd p;
\]
and in a similar way it will be
\[
1+q=(\surd p+\surd q)\surd q;
\]
with these values introduced, our sum turns into
\[
s=\frac{1}{2}(\surd p+\surd q)^n(p^{\frac{n}{2}}+q^{\frac{n}{2}}).
\]
Now, since in general it is
\[
p^\alpha+q^\alpha=2\Cos \alpha\varphi,
\]
it will be
\[
p^{\frac{1}{2}}+q^{\frac{1}{2}}=2\Cos\frac{1}{2}\varphi \quad \textrm{and}
\quad p^{\frac{n}{2}}+q^{\frac{n}{2}}=2\Cos \frac{1}{2}n\varphi,
\]
and with these values substituted, the sought for sum will now be expressed
by reals in the following way:
\[
s=2^n\Cos\frac{1}{2}\varphi^n\Cos \frac{1}{2}n\varphi.
\]

\S 6. With this understood, we obtain a most remarkable summation,
which is such that it is always
\[
\begin{split}
&1+\frac{n}{1}\Cos \varphi+\frac{n}{1}\cdot \frac{n-1}{2}\Cos 2\varphi
+\frac{n}{1}\cdot
\frac{n-1}{2}\cdot \frac{n-2}{3}\Cos 3\varphi+\textrm{etc.}\\
&=2^n\Cos \frac{1}{2}\varphi^n \Cos \frac{1}{2}n\varphi,
\end{split}
\]
the truth of which always holds for whatever number that is substituted
for $n$, either integers or fractions or even negatives. It will thus be worthwhile 
to relate to the eyes simpler cases from this type.

\begin{center}
{\Large Expansion of the cases in which the exponent $n$ is a positive integral number}
\end{center}

\S 7. Let us consider the following cases:

$1^\circ$. Let $n=0$, and the series becomes unity, and now the sum will be$=1$.

$2^\circ$. Let $n=1$, and the series becomes
\[
1+\Cos\varphi;
\]
and the sum is found to be
\[
2\Cos\frac{1}{2}\varphi^2,
\]
as it is apparent that
\[
2\Cos\frac{1}{2}\varphi^2=1+\Cos\varphi.
\]

$3^\circ$. Let $n=2$, and the series turns into
\[
1+2\Cos\varphi+\Cos 2\varphi,
\]
and the sum turns out to be
\[
=4\Cos \frac{1}{2}\varphi^2\Cos\varphi.
\]
For we have just seen that $2\Cos\frac{1}{2}\varphi^2=1+\Cos\varphi$,
and this form multiplied by $2\Cos\varphi$ produces
\[
2\Cos\varphi\cdot 2\Cos\frac{1}{2}\varphi^2=1+2\Cos\varphi+\Cos 2\varphi.
\]

$4^\circ$. Now let $n=3$, and this sum arises
\[
1+3\Cos\varphi+3\Cos 2\varphi+\Cos 3\varphi,
\]
whose sum is
\[
=8\Cos \frac{1}{2}\varphi^3\Cos \frac{3}{2}\varphi.
\]
This formula yields the above series by reductions well known enough.

$5^\circ$. Now let $n=4$, and the series turns into
\[
1+4\Cos\varphi+6\Cos 2\varphi+4\Cos 3\varphi+\Cos 4\varphi;
\]
whose sum will therefore be
\[
2^4\Cos \frac{1}{2}\varphi^4\Cos 2\varphi,
\]
whose truth can indeed be shown without difficulty. Thus one can always show the
truth by known reductions.

\begin{center}
{\Large Expansion of the cases in which a negative integral number is taken
for $n$}
\end{center}

\S 8. Let us first put $n=-1$, and the following infinite series arises
\[
1-\Cos\varphi+\Cos 2\varphi-\Cos 3\varphi+\Cos 4\varphi-5\Cos \varphi+\textrm{etc. to infinity},
\]
whose sum by our general series will therefore be
\[
\frac{\Cos \frac{1}{2}\varphi}{2\Cos\frac{1}{2}\varphi}=\frac{1}{2},
\]
which has previously been observed by Geometers. 
But if we let this series, whose sum we for the present put $=s$,
go to $2\Cos\frac{1}{2}\varphi$, it will be found by very
familiar reductions
\[
2s\Cos\frac{1}{2}\varphi=
\left\{ \begin{array}{l}
2\Cos\frac{1}{2}\varphi-\Cos\frac{3}{2}\varphi+\Cos\frac{5}{2}\varphi-\Cos\frac{7}{2}\varphi
+\Cos\frac{9}{2}\varphi-\\
-\Cos\frac{1}{2}\varphi+\Cos\frac{3}{2}\varphi-\Cos\frac{5}{2}\varphi
+\Cos\frac{7}{2}\varphi-\Cos\frac{9}{2}\varphi+
\end{array} \right\}
\textrm{etc.},
\]
which obviously reduces to $2s\Cos\frac{1}{2}\varphi=\Cos\frac{1}{2}\varphi$
and thus $s=\frac{1}{2}$.

\S 9. Let us now set $n=-2$, and the following series will arise
\[
1-2\Cos\varphi+3\Cos 2\varphi-4\Cos 3\varphi+5\Cos 4\varphi-6\Cos 5\varphi+\textrm{etc.},
\]
whose sum will therefore be
\[
=\frac{\Cos\varphi}{4\Cos\frac{1}{2}\varphi^2},
\]
the truth of which can even now be shown in the following way.
With the sum of the series put $=s$, it will be
\[
2s\Cos\frac{1}{2}\varphi=
\left\{ \begin{array}{l}
2\Cos\frac{1}{2}\varphi-2\Cos\frac{3}{2}\varphi+3\Cos\frac{5}{2}\varphi
-4\Cos\frac{7}{2}\varphi+\\
-2\Cos\frac{1}{2}\varphi+3\Cos\frac{3}{2}\varphi-4\Cos\frac{5}{2}\varphi
+5\Cos\frac{7}{2}\varphi-
\end{array} \right\} \textrm{etc.},
\]
whose value reduces to the following series
\[
2s\Cos\frac{1}{2}\varphi=\Cos\frac{3}{2}\varphi-\Cos\frac{5}{2}\varphi
+\Cos\frac{7}{2}\varphi-\Cos\frac{9}{2}\varphi+\textrm{etc.}
\]
Let us next multiply by $2\Cos\frac{1}{2}\varphi$, which yields
\[
4s\Cos\frac{1}{2}\varphi^2=\left\{
\begin{array}{l}
\Cos\varphi+\Cos 2\varphi-\Cos 3\varphi+\Cos 4\varphi-\Cos 5\varphi+\\
-\Cos 2\varphi+\Cos 3\varphi-\Cos 4\varphi+\Cos 5\varphi-
\end{array}\right\} \textrm{etc.}=\Cos\varphi
\]
and thus
\[
s=\frac{\Cos \varphi}{4\Cos\frac{1}{2}\varphi^2},
\]
as we had found; or, it is also
\[
s=\frac{\Cos\varphi}{2(1+\Cos\varphi)}.
\]

\S 10. Now let $n=-3$, and this infinite series will arise
\[
1-3\Cos\varphi+6\Cos 2\varphi-10\Cos 3\varphi+15\Cos 4\varphi-21\Cos 5\varphi+\textrm{etc.},
\]
whose sum will therefore be
\[
=\frac{\Cos\frac{3}{2}\varphi}{8\Cos\frac{1}{2}\varphi^3}.
\]
This expression however can in turn be reduced to
\[
s=\frac{1}{2}-\frac{3}{8\Cos\frac{1}{2}\varphi^2}=\frac{1}{2}-\frac{3}{4(1+\Cos\varphi)},
\]
and so that it is thus
\[
s=\frac{-1+2\Cos\varphi}{4(1+\Cos\varphi)}.
\]

\S 11. The following summations will be obtained in a similar way:
\begin{eqnarray*}
1-4\Cos\varphi+10\Cos 2\varphi-20\Cos 3\varphi+35\Cos 4\varphi-\textrm{etc.}&=&\frac{\Cos 2\varphi}{16\Cos \frac{1}{2}\varphi^4},\\
1-5\Cos\varphi+15\Cos 2\varphi-35\Cos 3\varphi+70\Cos 4\varphi-\textrm{etc.}&=&\frac{\Cos \frac{5}{2}\varphi}{32\Cos\frac{1}{2}\varphi^5},\\
1-6\Cos\varphi+21\Cos 2\varphi-56\Cos 3\varphi+126\Cos 4\varphi-\textrm{etc.}&=&\frac{\Cos 3\varphi}{64\Cos \frac{1}{2}\varphi^6},\\
1-7\Cos\varphi+28\Cos 2\varphi-84\Cos 3\varphi+210\Cos 4\varphi-\textrm{etc.}&=&\frac{\Cos\frac{7}{2}\varphi}{128\Cos \frac{1}{2}\varphi^7}\\
\textrm{etc.}&&
\end{eqnarray*}

\begin{center}
{\Large Expansion of the case in which $n=\frac{1}{2}$}
\end{center}

\S 12. Thus here the following infinite series will be formed
\[
1+\frac{1}{2}\Cos\varphi-\frac{1\cdot 1}{2\cdot 4}\Cos 2\varphi
+\frac{1\cdot 1\cdot 3}{2\cdot 4\cdot 6}\Cos 3\varphi-
\frac{1\cdot 1\cdot 3\cdot 5}{2\cdot 4\cdot 6\cdot 8}\Cos 4\varphi+\textrm{etc.},
\]
whose sum will therefore be
\[
=\Cos\frac{1}{4}\varphi\surd 2\Cos\frac{1}{2}\varphi,
\]
the truth of which will be not at all easy to justify otherwise;
however, certain cases come at once to the eyes. Truly if it were
$\varphi=0$, we obtain
\[
1+\frac{1}{2}-\frac{1\cdot 1}{2\cdot 4}+\frac{1\cdot 1\cdot 3}{2\cdot 4\cdot 6}
-\frac{1\cdot 1\cdot 3\cdot 5}{2\cdot 4\cdot 6\cdot 8}+\textrm{etc.}=\surd 2;
\]
which series clearly arises from the expansion
\[
(1+1)^\frac{1}{2}=\surd 2.
\]
Let us now make it $\varphi=180^\circ$, so that $\frac{1}{2}\varphi=90^\circ$,
and the series will be
\[
1-\frac{1}{2}-\frac{1\cdot 1}{2\cdot 4}-\frac{1\cdot 1\cdot 3}{2\cdot 4\cdot 6}
-\frac{1\cdot 1\cdot 3\cdot 5}{2\cdot 4\cdot 6\cdot 8}-\textrm{etc.}=0,
\]
which is also evident, because this series arises from the form $(1-1)^\frac{1}{2}$. Furthermore, let $\varphi=90^\circ$, and the series which then arises
will be
\[
1+\frac{1\cdot 1}{2\cdot 4}-\frac{1\cdot 1\cdot 3\cdot 5}{2\cdot 4\cdot 6\cdot 8}
+\frac{1\cdot 1\cdot 3\cdot 5\cdot 7\cdot 9}{2\cdot 4\cdot 6\cdot 8\cdot 10\cdot 12}-\textrm{etc.}
=\Cos 22^\circ 30' \sqrt[4]{2}.
\]
Indeed it is
\[
\Cos 22^\circ 30'=\surd \frac{(1+\Cos 45^\circ)}{2}=\surd \frac{1+\surd 2}{2\surd 2},
\]
from which the sum is concluded to be
\[
=\surd\frac{(1+\surd 2)}{2}
\]
and thus this very noteworthy summation is obtained
\[
1+\frac{1\cdot 1}{2\cdot 4}-\frac{1\cdot 1\cdot 3\cdot 5}{2\cdot 4\cdot 6\cdot 8}+
\frac{1\cdot 1\cdot 3\cdot 5\cdot 7\cdot 9}{2\cdot 4\cdot 6\cdot 8\cdot 10\cdot 12}-\textrm{etc.}=\surd \frac{1+\surd 2}{2}.
\]
Let us further put $\varphi=60^\circ$, and this series will arise
\[
\begin{split}
&1+\frac{1}{2}\cdot \frac{1}{2}+\frac{1\cdot 1}{2\cdot 4}\cdot \frac{1}{2}
-\frac{1\cdot 1\cdot 3}{2\cdot 4\cdot 6}\cdot 1+
\frac{1\cdot 1\cdot 3\cdot 5}{2\cdot 4\cdot 6\cdot 8}\cdot \frac{1}{2}
+\frac{1\cdot 1\cdot 3\cdot 5\cdot 7}{2\cdot 4\cdot 6\cdot 8\cdot 10}
\cdot \frac{1}{2}\\
&-\frac{1\cdot 1\cdot 3\cdot 5\cdot 7\cdot 9}{2\cdot 4\cdot 6\cdot 8\cdot 10\cdot 12}\cdot 1+\textrm{etc.},
\end{split}
\]
whose sum will therefore be $\Cos 15^\circ \sqrt[4]{3}$. Therefore since it is
\[
\Cos 15^\circ=\surd\frac{1+\Cos 30^\circ}{2}=\surd \frac{2+\surd 3}{4},
\]
the sum of the series will be
\[
\frac{1}{2}\surd(3+2\surd 3).
\]

\begin{center}
{\Large Expansion of the case in which $n=-\frac{1}{2}$}
\end{center}

\S 13. Thus here the following infinite series will be formed
\[
1-\frac{1}{2}\Cos \varphi+\frac{1\cdot 3}{2\cdot 4}\Cos 2\varphi
-\frac{1\cdot 3\cdot 5}{2\cdot 4\cdot 6}\Cos 3\varphi
+\frac{1\cdot 3\cdot 5\cdot 7}{2\cdot 4\cdot 6\cdot 8}\Cos 4\varphi-\textrm{etc.},
\]
whose sum will therefore be
\[
=\frac{\Cos \frac{1}{4}\varphi}{\surd(2\Cos\frac{1}{2}\varphi)}.
\]
So, if it were $\varphi=0$, this summation will follow
\[
1-\frac{1}{2}+\frac{1\cdot 3}{2\cdot 4}-\frac{1\cdot 3\cdot 5}{2\cdot 4\cdot 6}
+\frac{1\cdot 3\cdot 5\cdot 7}{2\cdot 4\cdot 6\cdot 8}-\textrm{etc.}=
\frac{1}{\surd 2}.
\]
For this series arises from the form $(1+1)^{-\frac{1}{2}}$. Now
let $\varphi=180^{\circ}$, and the resulting series will be
\[
1+\frac{1}{2}+\frac{1\cdot 3}{2\cdot 4}+\frac{1\cdot 3\cdot 5}{2\cdot 4\cdot 6}+
\frac{1\cdot 3\cdot 5\cdot 7}{2\cdot 4\cdot 6\cdot 8}+\textrm{etc.}=\infty.
\]
For this series follows from the expansion of $(1-1)^{-\frac{1}{2}}$.
Furthermore, let us take $\varphi=90^\circ$, and the series will be
\[
1-\frac{1\cdot 3}{2\cdot 4}+\frac{1\cdot 3\cdot 5\cdot 7}{2\cdot 4\cdot 6\cdot 8}
-\frac{1\cdot 3\cdot 5\cdot 7\cdot 9\cdot 11}{2\cdot 4\cdot 6\cdot 8\cdot 10\cdot 12}
+\textrm{etc.}=\frac{\Cos 22^\circ 30'}{\sqrt[4]{2}}.
\]
But we have seen before that
\[
\Cos 22^\circ 30'=\surd \frac{1+\surd 2}{2\surd 2},
\]
whence this sum will be
\[
=\frac{1}{2}\surd(1+\surd 2).
\]

\S 14. Also in general for all the exponents $n$, it will be worthwhile
to take certain values for the angle $\varphi$; and first indeed taking
$\varphi=0$ we will have
\[
1+\combin{n}{1}+\combin{n}{2}+\combin{n}{3}+\combin{n}{4}+\textrm{etc.}=2^n;
\]
this series is namely the expansion of the formula $(1+1)^n$.
Now let us take $\varphi=180^\circ$ and this series will arise
\[
1-\combin{n}{1}+\combin{n}{2}-\combin{n}{3}+\combin{n}{4}-\textrm{etc.}=0,
\]
namely this series is $=(1-1)^n$.
Further let $\varphi=90^\circ$, and then the arising series will be
\[
1-\combin{n}{2}+\combin{n}{4}-\combin{n}{6}+\combin{n}{8}-\combin{n}{10}+\textrm{etc.},
\]
whose sum will therefore be
\[
2^n(\Cos 45^\circ)^n\Cos n45^\circ =2^{\frac{1}{2}n}\Cos n45^\circ.
\]

\S 15. This last series seems worth more attention, because the truth of it
is in no small way concealed; 
a few special cases from it will be considered, and first indeed for
positive integral numbers.

\begin{tabular}{ll}
$1^\circ$.&If $n=0$, it will be $1=1$,\\
$2^\circ$.&If $n=1$, it will be $1=\Cos 45^\circ \surd 2=1$,\\
$3^\circ$.&If $n=2$, it will be $1-1=2\Cos 90^\circ=0$,\\
$4^\circ$.&If $n=3$, it will be $1-3=2^\frac{3}{2}\Cos 3\cdot 45^\circ=-2$,\\
$5^\circ$.&If $n=4$, it will be $1-6+1=2^2\Cos 4\cdot 45^\circ=-4$,\\
$6^\circ$.&If $n=5$, it will be $1-10+5=2^\frac{5}{2}\Cos 5\cdot 45^\circ=-4$,\\
$7^\circ$.&If $n=6$, it will be $1-15+15-1=2^3\Cos 6\cdot 45^\circ=0$,\\
$8^\circ$.&If $n=7$, it will be $1-21+35-7=2^\frac{7}{2}\Cos 7\cdot 45^\circ=2^3$,\\
$9^\circ$.&If $n=8$, it will be $1-28+70-28+1=2^4\Cos 8\cdot 45^\circ=2^4$\\
&\textrm{etc.}
\end{tabular}

\S 16. The cases merit more attention in which negative numbers
are taken for $n$,
obviously in which infinite series occur.

$1^\circ$. If $n=-1$, it will be
\[
1-1+1-1+1-1+1-1+1-1+\textrm{etc.}=\frac{\Cos 45^\circ}{\surd 2}=\frac{1}{2},
\]

$2^\circ$. If $n=-2$, it will be
\[
1-3+5-7+9-11+13-15+17-\textrm{etc.}=\frac{\Cos 2\cdot 45^\circ}{2}=0,
\]

$3^\circ$. If $n=-3$, it will be
\[
1-6+15-28+45-66+91-\textrm{etc.}=\frac{\Cos 3\cdot 45^\circ}{\surd 8}=-\frac{1}{4},
\]

$4^\circ$. If $n=-4$, it will be
\[
1-10+35-84+165-286+455-\textrm{etc.}=\frac{\Cos 4\cdot 45^\circ}{4}=-\frac{1}{4},
\]

$5^\circ$. If $n=-5$, it will be
\[
1-15+70-210+495-1001+\textrm{etc.}=\frac{\Cos 5\cdot 45^\circ}{2^\frac{5}{2}}=-\frac{1}{8},
\]

$6^\circ$. If $n=-6$, it will be
\[
1-21+126-462+1287-3003+\textrm{etc.}=\frac{\Cos 6\cdot 45^\circ}{8}=0
\]

\begin{center}
etc.
\end{center}

For these cases it will be in general
\[
\begin{split}
&1-\frac{\lambda(\lambda+1)}{1\cdot 2}+\frac{\lambda(\lambda+1)(\lambda+2)(\lambda+3)}{1\cdot 2\cdot 3\cdot 4}
-\frac{\lambda(\lambda+1)(\lambda+2)(\lambda+3)(\lambda+4)(\lambda+5)}{1\cdot
2\cdot 3\cdot 4\cdot 5}\\
&+\frac{\lambda(\lambda+1)\cdots(\lambda+7)}{1\cdot 2\cdots 8}
-\frac{\lambda(\lambda+1)\cdots (\lambda+9)}{1\cdot 2\cdots 10}+\textrm{etc.},
\end{split}
\]
thus the sum of this series will be
\[
=\frac{\Cos \lambda\cdot 45^\circ}{2^{\frac{1}{2}\lambda}}.
\]

\begin{center}
{\Large Problem 2}
\end{center}

{\em Given this series of sines
\[
\frac{n}{1}\Sin \varphi+\frac{n}{1}\cdot \frac{n-1}{2}\Sin 2\varphi
+\frac{n}{1}\cdot \frac{n-1}{2}\cdot \frac{n-2}{3}\Sin 3\varphi+\textrm{etc.}=s,
\]
so that with the characters introduced above it will be
\[
s=\combin{n}{0}\Sin 0\varphi+\combin{n}{1}\Sin \varphi+\combin{n}{2}\Sin 2\varphi
+\combin{n}{3}\Sin 3\varphi+\textrm{etc.},
\]
to express the value this sum $s$ with reals.}

\begin{center}
{\Large Solution}
\end{center}

\S 17. If we introduce here again the letters 
\[
p=\Cos \varphi+\surd -1\cdot \Sin \varphi
\]
and
\[
q=\Cos \varphi-\surd -1\cdot \Sin \varphi,
\]
since it is
\[
p^n-q^n=2\surd -1\cdot \Sin n\varphi,
\]
the given series can be separated into the following two
\[
2s\surd -1=
\left\{ \begin{array}{l}
+\combin{n}{1}p+\combin{n}{2}pp+\combin{n}{3}p^3+\combin{n}{4}p^4+\\
-\combin{n}{1}q-\combin{n}{2}qq-\combin{n}{3}q^3-\combin{n}{4}q^4-
\end{array}\right\}
\textrm{etc.},
\]
from which clearly it will be
\[
2s\surd -1=(1+p)^n-(1+q)^n.
\]

\S 18. It will be helpful to again observe now that
\[
1+p=(\surd p+\surd q)\surd p
\]
and
\[
1+q=(\surd p+\surd q)\surd q.
\]
Employing these values it will be
\[
2s\surd -1=(\surd p+\surd q)^n(p^{\frac{n}{2}}-q^{\frac{n}{2}}).
\]
Thus, since it is
\[
p^{\frac{n}{2}}-q^{\frac{n}{2}}=2\surd -1\cdot \Sin \frac{1}{2}n\varphi
\]
and
\[
\surd p+\surd q=2\Cos \frac{1}{2}\varphi,
\]
dividing by $2\surd -1$ yields a real expression for the sought for sum
\[
s=2^n\Cos \frac{1}{2}\varphi^n \Sin \frac{1}{2}n\varphi.
\]

\end{document}